\documentclass[10pt]{article}
\usepackage{amsfonts,amssymb,amsbsy,latexsym,amsthm,amsmath,tabulary,graphicx,times,caption,fancyhdr}
\usepackage[utf8]{inputenc}
\usepackage{amscd, graphicx, color}
\usepackage{cite}
\usepackage{url,multirow,morefloats,floatflt,cancel,tfrupee}

\usepackage{enumerate}
\usepackage{enumitem}
\usepackage{authblk}
\usepackage{cases}
\usepackage[bookmarksnumbered, plainpages]{hyperref}
\usepackage{fancyhdr}
\usepackage{caption}
\usepackage{subcaption}
\usepackage{cite}
\usepackage[a4paper,bindingoffset=0.2in,%
left=1in,right=1in,top=1in,bottom=1in,%
footskip=0.30in]{geometry}
\usepackage{etoolbox}

\usepackage{xargs} 
\usepackage{bbm} 
\usepackage{stackrel} 
\usepackage{xspace}  
\usepackage{diagbox} 
\usepackage{calligra}
\usepackage{calrsfs}
\usepackage[mathscr]{euscript}
\newtheorem{theorem}{Theorem}[section]

\newtheorem{proposition}[theorem]{Proposition}
\newtheorem{notation}[theorem]{Notation}

\theoremstyle{Definition}

\newtheorem{definition}[theorem]{Definition}

\usepackage{xargs} 
\usepackage{bbm} 
\usepackage{stackrel} 
\usepackage{xspace}  
\usepackage{diagbox} 
\usepackage{float}  
\usepackage{xpatch}
\makeatletter
\usepackage{framed}  
\usepackage{xcolor}    
\usepackage{colortbl}  
\usepackage{listings} 
\definecolor{grigiochiarissimo}{rgb}{0.97,0.97,0.97}
\definecolor{grigiomoltochiaro}{rgb}{0.96,0.96,0.96}
\definecolor{grigiochiaro}{rgb}{0.8,0.8,0.8}
\definecolor{grigiomedio}{rgb}{0.6,0.6,0.6}
\definecolor{grigioscuro}{rgb}{0.3,0.3,0.3}
\definecolor{giallochiaro}{rgb}{1,1,0.95}
\definecolor{azzurrochiaro}[h=20]{rgb}{0.95,0.97,1}
\definecolor{blu}{rgb}{0,0,0.8}
\definecolor{bluscuro}{rgb}{0,0,0.5}
\definecolor{verdescuro}{rgb}{0,0.5,0}
\definecolor{arancione}{rgb}{1,0.5,0}
\definecolor{viola}{rgb}{0.65, 0.12, 0.82}
\definecolor{malva}{rgb}{0.58,0,0.82}
\definecolor{coloredisfondo}{rgb}{0.97,0.97,0.97}         
\definecolor{colorebordo}{HTML}{E0E0E0}                 
\definecolor{coloredibase}{HTML}{000000}                
\definecolor{colorenumerazione}{rgb}{0.8,0.8,0.8}       
\definecolor{coloreparolechiave}{rgb}{0,0,0.8}          
\definecolor{coloreparolechiave-2}{rgb}{0.65,0.12,0.82} 
\definecolor{coloreparolechiave-3}{HTML}{002DB2}        
\definecolor{coloreparolechiave-4}{HTML}{00238C}        
\definecolor{coloreparolechiave-5}{HTML}{00238C}        
\definecolor{colorestringhe}{rgb}{0,0.5,0}              
\definecolor{coloreidentificatori}{HTML}{000000}        
\definecolor{colorestringhe}{HTML}{067D17}              
\definecolor{coloredirettive}{HTML}{9D8740}             
\definecolor{colorenumero}{HTML}{001A66}                
\definecolor{coloredelimitatori}{HTML}{FF0000}
\definecolor{coloredelimitatoriangolari}{rgb}{0.65, 0.12, 0.82}  
\definecolor{coloredelimitatoriquadre}{rgb}{0.65, 0.12, 0.82}    
\definecolor{colorecommenti}{rgb}{0.6,0.6,0.6}                   
\definecolor{colorecommenti-2}{HTML}{4497E7}                     
\lstdefinestyle{stilePython}{
language=Python,%
keywordstyle={\bfseries\color{coloreparolechiave}},            
morekeywords={},          
keywords=[2]{self,len},                                        
keywordstyle={[2]{\color{coloreparolechiave-2}}},              
keywords=[3]{},                                                
keywordstyle={[3]{\color{coloreparolechiave-3}}},              
keywords=[4]{},                                                
keywordstyle={[4]{\color{coloreparolechiave-4}}},              
keywords=[5]{},                                                
keywordstyle={[5]{\color{coloreparolechiave-5}}},              
identifierstyle=\color{coloreidentificatori},                  
stringstyle={\ttfamily \color{colorestringhe}},                
commentstyle=\color{colorecommenti},                           
morecomment=[s][\color{colorecommenti-2}]{/*}{*/},             
literate=%
{0}{{\color{colorenumero}{0}}}{1}%
{1}{{\color{colorenumero}{1}}}{1}%
{2}{{\color{colorenumero}{2}}}{1}%
{3}{{\color{colorenumero}{3}}}{1}%
{4}{{\color{colorenumero}{4}}}{1}%
{5}{{\color{colorenumero}{5}}}{1}%
{6}{{\color{colorenumero}{6}}}{1}%
{7}{{\color{colorenumero}{7}}}{1}%
{8}{{\color{colorenumero}{8}}}{1}%
{9}{{\color{colorenumero}{9}}}{1},                 
backgroundcolor=\color{coloredisfondo},                        
frame=none,                                                      
xleftmargin=8pt,  
xrightmargin=4pt,  
framesep=2pt,                                                  
framerule=0.6pt,                                                 
rulecolor=\color{colorebordo},                                 
numbers=left,                                                  
numberstyle=\tiny\color{colorenumerazione},                    
numbersep=3pt,                                                 
numberblanklines=false,                                        
showlines=false,                                               
firstnumber=1,                                                 
breakatwhitespace=false,                                       
breaklines=true,                                               
captionpos=b,                                                  
keepspaces=true,                                               
showspaces=false,                                              
showstringspaces=false,                                        
showtabs=false,                                                
columns=fixed,                                                 
tabsize=3,                                                     
extendedchars=true,                                            
escapeinside=$$                     
}
\lstdefinestyle{stilePythonambiente}{
style=stilePython,
basicstyle={\footnotesize\ttfamily\color{coloredibase}},
firstnumber=1,
lineskip=-1pt
}
\lstdefinestyle{stilesemplicecomando}{
basicstyle={\ttfamily\color{coloredibase}},
mathescape
}
\lstdefinestyle{stilePythoncomando}{
style=stilePython,
basicstyle={\ttfamily\color{coloredibase}},
mathescape
}
\lstdefinestyle{stilePythonconsole}{
style=stilePython,
numbers=none,
backgroundcolor=\color{giallochiaro},
basicstyle={\footnotesize\ttfamily\color{coloredibase}},
lineskip=-2pt,
otherkeywords={>>>}
}
\newcommand{\cpy}[1]{\lstinline[style=stilePythoncomando]!#1!}    

\lstnewenvironment{codpy}{%
\lstset{style=stilePythonambiente}}{}
\lstnewenvironment{codpyconsole}{%
\lstset{style=stilePythonconsole}}{}

\usepackage{pgf-umlcd}   

\tikzset{
  reduce height/.style={
    minimum height=0pt,
    inner ysep=0pt,
    text depth=2pt
  },
  reduce height/.default={0pt}
}

\newcommand{\df}[1]{\textit{#1}}  
\newcommand{\Tau}{\mathcal{T}} 
\newcommand{\U}{\mathbbmss{U}}

\newcommand{\NN}{\mathbb{N}} 

\newcommand{\cS}{\mathcal{S}}
\newcommand{\cA}{\mathcal{A}}
\newcommand{\cB}{\mathcal{B}}

\newcommand{\und}{\underline}  

\newcommand{\M}[1]{\mu_{#1}}       
\newcommand{\I}[1]{\sigma_{#1}}       
\newcommand{\NM}[1]{\omega_{#1}}   
\newcommand{\DM}[2][u]{\mu_{#2}\left(#1\right)}       
\newcommand{\DI}[2][u]{\sigma_{#2}\left(#1\right)}       
\newcommand{\DNM}[2][u]{\omega_{#2}\left(#1\right)}   
\newcommand{\ANG}[1]{\left\langle #1 \right\rangle}  
\newcommand{\NSbase}[4][\U]{\ANG{#1, #2, #3, #4}}  
\newcommand{\ns}[1]{\widetilde{#1}}   
\newcommand{\NS}[2][\U]{\NSbase[#1]{\M{#2}}{\I{#2}}{\NM{#2}}}   
\newcommand{\nNS}[2][\U]{\ns{#2} = \NS[#1]{#2}}  
\newcommandx{\NSEXT}[3][1=u,2=\U]{\left\{ \left( #1, \DM{#3}, \DI{#3}, \DNM{#3} \right) : \, #1 \in #2 \right\}} 
\newcommand{\nameSVNS}[1][]{SVN-set#1\xspace}  
\newcommand{\NSemptyset}[1][]{\widetilde{\emptyset}_{#1}}  
\newcommand{\NSabsoluteset}[1][\U]{\widetilde{#1}}         
\newcommand{\NSsubseteq}{\mathrel{\ooalign{
\raise0.2ex\hbox{$\Subset$}%
\cr\hidewidth%
\raise-0.25ex\hbox{\rule[0.8pt]{4pt}{0.4pt}}%
\hidewidth\cr%
}}}
\newcommand{\NSsupseteq}{\mathrel{\ooalign{
\raise0.2ex\hbox{$\Supset$}%
\cr\hidewidth%
\raise-0.25ex\hbox{\rule[0.8pt]{4pt}{0.4pt}}%
\hidewidth\cr%
}}}
\newcommand{\NSeq}{\mathrel{\ooalign{
\raise0.1ex\hbox{$=$}%
\cr%
\hskip0.6pt\raise0.2ex\hbox{\rule{4.2pt}{0.35pt}}
\cr\hidewidth%
\raise1.1ex\hbox{{\rule{4.2pt}{0.35pt}}}
\hidewidth\cr%
}}}
\newcommand{\NScompl}{{\hskip.15ex\ooalign{\hbox{\scalebox{0.8}{$\complement$}}%
\hidewidth\cr\hspace{1.3pt}\hbox{\rule[-0.2pt]{0.5pt}{5.0pt}}%
}}}    
\newcommand{\NScup}{\ensuremath{\Cup}}  
\newcommand{\NScap}{\ensuremath{\Cap}}  
\newcommand{\NSCup}{
\mathop{\vphantom{\bigcup}\vcenter{\hbox{\text{%
\ooalign{$\displaystyle\bigcup$\cr%
\hidewidth%
\raisebox{.25ex}{\resizebox{.7\width}{.875\height}{$\displaystyle\bigcup$}}%
\hidewidth\cr
}}}}}}%
\newcommand{\NSCap}{
\mathop{\vphantom{\bigcap}\vcenter{\hbox{\text{%
\ooalign{$\displaystyle\bigcap$\cr%
\hidewidth%
\raisebox{-.1ex}{\resizebox{.7\width}{.875\height}{$\displaystyle\bigcap$}}%
\hidewidth\cr%
}}}}}}%

\newcommand{\dunder}{\ensuremath{\_\hspace{0.2pt}\_}}   
\newcommand{\nsempty}{\ensuremath{\widetilde{\emptyset}}}   
\newcommand{\nsabsol}{\ensuremath{\widetilde{\mathbb{U}}}}   
\newcommand{\tld}{\ensuremath{\sim}}  
\newcommand{\ds}{\displaystyle}

\newenvironment{enumi}[1][0]%
{\begin{enumerate}[label={\rm(\arabic*)},topsep=4pt]%
\itemsep1pt
\setcounter{enumi}{#1}%
}%
{\end{enumerate}}

{\begin{enumerate}[label={\rm(\alph*)},topsep=4pt]%
\itemsep1pt
\setcounter{enumi}{#1}%
}%
{\end{enumerate}}

{\begin{enumerate}[label={\rm(\roman*)},topsep=4pt]%
\itemsep1pt
\setcounter{enumi}{#1}%
}%
{\end{enumerate}}

{\begin{tabular}[t]{|c|c|c|c|}\hline
\diagbox{\quad$#1$}{\\[1mm]\quad$#2$} &
$#3$ & $#4$ & $#5$ \\ \hline}%
{\\ \hline \end{tabular}}

\usepackage{dirtree}   
\DTsetlength{0.2em}{20pt}{0.2em}{0.8pt}{4pt}

\usepackage{setspace}  

\usepackage[ruled,vlined]{algorithm2e}  
\SetAlgoNlRelativeSize{-1}
\newenvironment{algoritmo}{%
\begin{algorithm}[H]%
\setstretch{0.9}%
\NoCaptionOfAlgo  
\LinesNotNumbered 
\DontPrintSemicolon 
}
{\end{algorithm}}

\begin{document}
\setcounter{page}{08}
\begin{center}
\textbf{\Large A Python Framework Enhancement for Neutrosophic Topologies}\\
\vspace{10pt}
\textbf{Giorgio Nordo$^{1}$, Saeid Jafari$^{2}$ and Maikel Yelandi Leyva Vázquez$^{3}$}\\
\vspace{10pt}
$^{1}$ MIFT Department, University of Messina - Viale F. Stagno d'Alcontres 31, 98166 Messina, Italy.\\
e-mail: giorgio.nordo@unime.it
\\[2mm]
$^{2}$  Mathematical and Physical Science Foundation, Sidevej 5, 4200 Slagelse, Denmark. \\
e-mail: jafaripersia@gmail.com; saeidjafari@topositus.com
\\[2mm]
$^{3}$  Salesian Polytechnic University of Guayaquil, Guyas, Ecuador.\\
e-mail: mleyvaz@gmail.com\\
\end{center}

\section*{Abstract}
This paper introduces an extension to the Python Neutrosophic Sets (PYNS) framework,
originally detailed in \cite{nordo2024}, with the addition of the \texttt{NSfamily} class
for constructing and manipulating neutrosophic topologies. Building on existing classes
like \texttt{NSuniverse} and \texttt{NSset}, the \texttt{NSfamily} class enables
the definition and testing of neutrosophic families as basis and sub-basis for neutrosophic topological spaces.
This extension provides tools for verifying closure properties under union and intersection,
and for determining whether a given family constitutes a neutrosophic topology.
Through implemented algorithms, the framework automates the generation of topologies
from families of neutrosophic sets, offering an efficient tool for advancing research
in neutrosophic topology. Practical applications are demonstrated with detailed examples,
showcasing how this class enhances the scope and flexibility of neutrosophic modeling within the PYNS framework.
\\[4mm]
\noindent\textbf{Keywords:} Neutrosophic topology, Neutrosophic set, Python framework, Neutrosophic family, Neutrosophic topological basis, Neutrosophic topological sub-basis.
\\[4mm]
\noindent\textbf{2020 Mathematics Subject Classification:} 03E72, 54A40, 54D99, 68W99.


\section{Introduction}

The theory of neutrosophic sets, introduced by Smarandache in 1999 \cite{smarandache}, generalizes both fuzzy sets \cite{zadeh} and Atanassov's intuitionistic fuzzy sets \cite{atanassov} by providing three distinct parameters—truth, indeterminacy, and falsity—for each element. This generalization offers a more flexible and nuanced model for representing uncertainty and incomplete information, making it applicable to a wide range of fields, including statistics \cite{smarandache2014} and image processing \cite{zhang}. In addition, single-valued neutrosophic sets have been explored to represent degrees of truth, indeterminacy, and falsity within a simplified framework \cite{wang}, with applications in graph theory and decision making \cite{broumi, mondal}.

In recent years, neutrosophic topology has emerged as an extension of classical topological concepts to the neutrosophic domain. Key developments include the exploration of connectedness and stratification in neutrosophic topological spaces by Broumi et al. \cite{broumi2020} and the investigation of separation axioms by Dey and Ray \cite{dey2023}.
Additionally, foundational studies by Gallego Lupiáñez \cite{gallego2008, gallego2017,ray2021}
have laid the groundwork for applying neutrosophic theory to topology, creating opportunities to model
complex relationships characterized by uncertainty.

The practical application of neutrosophic sets has been greatly enhanced by the development of computational tools. El-Ghareeb introduced an open-source Python package to handle neutrosophic data \cite{el-ghareeb2019}, and Sleem et al. developed the PyIVNS tool for interval-valued neutrosophic operations \cite{sleem2020}. In addition, Topal et al. created a Python tool for implementing operations on bipolar neutrosophic matrices, useful in applications requiring complex matrix computations \cite{topal2019}. Furthermore, Saranya et al. developed a C\# application specifically designed to handle neutrosophic $g\alpha$-closed sets, expanding the accessibility of neutrosophic topology across different programming environments \cite{saranya2020}.
However, these tools are limited in scope and focus primarily on specific neutrosophic operations rather than providing
a comprehensive and interactive framework for dealing with symbolic representations of neutrosophic sets.

The Python Neutrosophic Sets (PYNS) framework, developed by Nordo et al. \cite{nordo2024}, addresses this limitation by offering a robust, modular library for the representation and manipulation of neutrosophic sets and mappings.
PYNS enables the modeling of various neutrosophic structures, such as single-valued neutrosophic mappings \cite{latreche2020} and neutrosophic soft topological spaces \cite{mehmood2020, mehmood2020a}. However, while PYNS supports essential topological operations, such as basis and sub-basis construction, it lacks a generalized framework for systematically defining and verifying neutrosophic topologies.

In this paper, we introduce the \texttt{NSfamily} class as an extension of the PYNS framework, aimed at facilitating the construction, manipulation, and verification of neutrosophic topologies. This new class draws upon theoretical work by Ozturk \cite{ozturk2020} and Salama et al. \cite{salama2012, salama2014, salama2014b, salama2014c} on neutrosophic sets and topological structures, offering a more integrated approach for defining families of neutrosophic sets as basis and sub-basis. The \texttt{NSfamily} class supports automated generation of neutrosophic topologies, including operations to verify closure under union and intersection, providing an efficient and scalable solution for neutrosophic topological modeling.

Furthermore, \texttt{NSfamily} addresses a broad range of computational requirements in neutrosophic topology, including the ability to handle interval-valued data and complex neutrosophic structures. This flexibility makes it suitable for both theoretical research and practical applications, as demonstrated in Rabuni and Balaman's approach to neutrosophic soft topologies using Python \cite{rabuni2023}, and in applications of single-valued neutrosophic ideals by Saber et al. \cite{saber2020}. Our contribution extends the PYNS framework with a cohesive and efficient tool for neutrosophic topological modeling, addressing the computational needs of researchers and expanding the possibilities of neutrosophic analysis in complex domains.

The structure of the paper is as follows. In Section 2, we review fundamental definitions and properties in neutrosophic topology, including closure and continuity. Section 3 details the \texttt{NSfamily} class, highlighting the core methods and their applications. Finally, we provide practical examples that illustrate the enhanced capabilities of PYNS with the \texttt{NSfamily} extension, discussing the potential impact of this work on future research in neutrosophic topology.


\section{Preliminaries}

This section provides a comprehensive overview of neutrosophic sets, neutrosophic topology, and the structure of the Python Neutrosophic Sets (PYNS) framework. This background is necessary to understand the design and functionality of the new \texttt{NSfamily} class, which extends PYNS for advanced neutrosophic topological modeling.

\subsection{Neutrosophic Sets}

The concept of neutrosophic sets, introduced by Smarandache \cite{smarandache}, generalizes classical, fuzzy \cite{zadeh}, and intuitionistic fuzzy sets \cite{atanassov} by assigning each element three independent degrees: truth, indeterminacy, and non-membership, each within the hyperreal interval $]0^{-},1^{+}[$ of the nonstandard real numbers.
A simpler variant, the single-valued neutrosophic set \cite{wang}, uses the standard interval
$[0,1]$, making it more accessible and practical for scientific and engineering applications.

\begin{definition}{\rm \cite{wang}}
\label{def:singlevaluedneutrosophicset}
Let \( \U \) be a universal set and \( A \subseteq \U \). A \df{single-valued neutrosophic set} (abbreviated \nameSVNS) over \( \U \), denoted by \( \nNS{A} \), is defined as:
\[
\ns{A} = \NSEXT{A}
\]
where \(\M{A} : \U \to I\), \(\I{A} : \U \to I\), and \(\NM{A} : \U \to I\) are the \df{membership}, \df{indeterminacy}, and \df{non-membership} functions of \(A\), respectively, and \(I = [0,1]\) denotes the real unit interval.
For every \( u \in \U \), \(\DM{A}\), \(\DI{A}\), and \(\DNM{A}\) are called the \df{degree of membership}, \df{degree of indeterminacy}, and \df{degree of non-membership} of \( u \), respectively.
\end{definition}

\begin{definition}{\rm\cite{smarandache,wang}}
\label{def:neutrosophicsubset_and_neutrosophicequality}
Let $\nNS{A}$ and $\nNS{B}$ be two \nameSVNS[s] over the universe set $\U$,
we say that:
\begin{itemize}
\item $\ns{A}$ is a \df{neutrosophic subset} (or simply a subset) of $\ns{B}$
and we write $\ns{A} \NSsubseteq \ns{B}$
if, for every $u \in \U$, it results
$\DM{A} \le \DM{B}$,
$\DI{A} \le \DI{B}$
and $\DNM{A} \ge \DNM{B}$.
We also say that $\ns{A}$ is contained in $\ns{B}$ or that $\ns{B}$ contains $\ns{A}$
and we write $\ns{B} \NSsupseteq \ns{A}$ to denote that $\ns{B}$ is a \df{neutrosophic superset} of $\ns{A}$.
\item $\ns{A}$ is a \df{neutrosophically equal} (or simply equal) to $\ns{B}$ and we write
$\ns{A} \NSeq \ns{B}$
if $\ns{A} \NSsubseteq \ns{B}$ and $\ns{B} \NSsubseteq \ns{A}$.
\end{itemize}
\end{definition}


\begin{notation}
Let $\U$ be a set, $I=[0,1]$ the unit interval of the real numbers, for every $r \in I$,
with $\und{r}$ we denote the constant mapping $\und{r}: \U \to I$
defined by $\und{r}(u) = r$, for every $u \in \U$.
\end{notation}

\begin{definition}{\rm\cite{wang}}
\label{def:neutrosophicemptyset_and_neutrosophicabsoluteset}
Given a universe set $\U$:
\begin{itemize}
\item the \nameSVNS $\NSbase{\und{0}}{\und{0}}{\und{1}}$
is said to be the \df{neutrosophic empty set} over $\U$
and it is denoted by $\NSemptyset$,
or more precisely by $\NSemptyset[\U]$ in case it is necessary
to specify the corresponding universe set
\item the \nameSVNS $\NSbase{\und{1}}{\und{1}}{\und{0}}$
is said to be the \df{neutrosophic absolute set} over $\U$
and it is denoted by $\NSabsoluteset$.
\end{itemize}
\end{definition}


On neutrosophic sets, union, intersection and complement operations can be defined.

\begin{definition}{\rm\cite{salama2013}}
\label{def:neutrosophic_union_and_inttersection}
Let $\left\{ \ns{A}_\alpha \right\}_{\alpha\in \Lambda}$ be a family of \nameSVNS[s]
$\ns{A}_\alpha=\NS{A_\alpha}$ over a common universe set $\U$, then:
\begin{itemize}
\item the \df{neutrosophic union}, denoted by $\ds \NSCup_{\alpha\in \Lambda} \ns{A}_\alpha$,
is the neutrosophic set $\nNS{A}$ with
$\ds \M{A} = \bigvee_{\alpha\in \Lambda} \M{A_\alpha}$,
$\ds \I{A} = \bigvee_{\alpha\in \Lambda} \I{A_\alpha}$, and
$\ds \NM{A} = \bigwedge_{\alpha\in \Lambda} \NM{A_\alpha}$.
In particular, the neutrosophic union of two \nameSVNS[s]
$\nNS{A}$ and $\nNS{B}$, denoted by $\ns{A} \NScup \ns{B}$, is the neutrosophic set defined by
$\NSbase{\M{A} \vee \M{B}}{\I{A} \vee \I{B}}{\NM{A} \wedge \NM{B}}$
\item the \df{neutrosophic intersection}, denoted by $\ds \NSCap_{\alpha\in \Lambda} \ns{A}_\alpha$,
is the neutrosophic set $\nNS{A}$ with
$\ds \M{A} = \bigwedge_{\alpha\in \Lambda} \M{A_\alpha}$,
$\ds \I{A} = \bigwedge_{\alpha\in \Lambda} \I{A_\alpha}$, and
$\ds \NM{A} = \bigvee_{\alpha\in \Lambda} \NM{A_\alpha}$.
In particular, the neutrosophic intersection of two \nameSVNS[s]
$\nNS{A}$ and $\nNS{B}$, denoted by $\ns{A} \NScap \ns{B}$, is the neutrosophic set defined by
$\NSbase{\M{A} \wedge \M{B}}{\I{A} \wedge \I{B}}{\NM{A} \vee \NM{B}}$
\end{itemize}
where $\wedge$ and $\vee$ denote the minimum and the maximum, respectively.
\end{definition}

\begin{definition}{\rm \cite{smarandache, wang}}
\label{def:neutrosophiccomplement}
The \df{neutrosophic complement} of a \nameSVNS \( \nNS{A} \) over \( \U \), denoted by \( \ns{A}^\NScompl \), is given by:
$
\ns{A}^\NScompl = \NSbase{\NM{A}}{\und{1} - \I{A}}{\M{A}}.
$
\end{definition}



\subsection{Neutrosophic Topology}

Neutrosophic topology generalizes classical topology by defining open and closed sets through neutrosophic set theory,
enabling topologies where openness, closedness, and indeterminacy are graded.
According to Salama et al \cite{salama2012}, Serkan Karatas et al \cite{karatas2016} and Gallego Lupiáñez \cite{gallego2008, gallego2017},
a neutrosophic topology is a family of neutrosophic sets which contains the empty neutrosophic set, the absolute neutrosophic set
and is closed with respect to the neutrosophic union and the finite neutrosophic intersection.

\begin{definition}\cite{gallego2008,salama2012}
\label{def:neutrosophictopology}
Let $\Tau$ be a family of neutrosophic sets over the same universe set $\U$, we say that
$\Tau$ is a \df{neutrosophic topology} if the following four conditions hold:
\begin{enumi}
\item $\NSemptyset \in \Tau$
\item $\NSabsoluteset \in \Tau$
\item $\forall \cA \subseteq \Tau$, $\NSCup \cA \in \Tau$
\item $\forall \ns{U}, \ns{V} \in \Tau$, $\ns{U} \NScap \ns{V} \in \Tau$
\end{enumi}
and when this occurs we also say that the pair $\left( \U, \Tau \right)$
constitutes a \df{neutrosophic topological space},
while the elements of $\Tau$ are named \df{neutrosophic open sets}.
\end{definition}

Although, the authors are not aware of papers in which notions such as bases and sub-basis
of neutrosophic topologies are introduced and treated, by simple analogy with classical topological spaces,
the following general definitions and basic properties can be introduced.

\begin{definition}
\label{def:comparisonofneutrosophictopologies}
Given two nutrosophic topologies $\Tau_1$ and $\tau_2$ over same universe set $\U$,
we say that $\Tau_1$ is \df{coarser} (weaker or smaller) than $\Tau_2$,
or equivalently, that $\Tau_2$ is \df{finer} (stronger or larger) than $\Tau_1$
if $\Tau_1 \subseteq \Tau_2$.
\end{definition}

So, the inclusion relation $\subseteq$ is a partial order on the set of all neutrosophic topologies over a universe set $\U$.

\begin{proposition}
\label{pro:intersectionofneutrosophictopologies}
The intersection $\ds \bigcap_{\alpha \in \Lambda} \Tau_\alpha$
of any family $\{\Tau_\alpha\}_{\alpha \in \Lambda}$ of neutrosophic topologies
over a common universe set $\U$ is itself a neutrosophic topology on $\U$.
\end{proposition}

\begin{proposition}
\label{pro:neutrosophictopologygeneratedbyafamily}
Given a family $\cS$ of neutrosophic subsets over a common universe set $\U$,
there exists and is unique the smallest (minimal) neutrosophic topology $\Tau(\cS)$ on $\U$ containing $\cS$.
\end{proposition}

\begin{definition}
\label{def:neutrosophictopologygeneratedbyafamily}
Let $\cS$  be a family of neutrosophic subsets over a common universe set $\U$,
the neutrosophic topology $\Tau(\cS)$ of Proposition \ref{pro:neutrosophictopologygeneratedbyafamily}
is called the \df{neutrosophic topology generated} by the family $\cS$.
\end{definition}

\begin{proposition}
\label{pro:neutrosophic_topologies_generated_by_two_families}
If $\cS_1$ and $\cS_2$ are two families of neutrosophic sets on the same universe set $\U$
such that $\cS_1 \subseteq \cS_2$
then the corresponding generated neutrosophic topologies satisfy $\Tau(\cS_1) \subseteq \Tau(\cS_2)$.
\end{proposition}

\begin{definition}
\label{def:neutrosophicbase}
Given a neutrosophic topological space $(\U, \Tau)$, a family $\cB \subseteq \Tau$
is said to be a \df{neutrosophic basis} for the neutrosophic topology $\Tau$
if every neutrosophic open set in $\Tau$ can be expressed as a neutrosophic union of a subfamily of $\cB$,
i.e., if for each $\ns{U} \in \Tau$, there exists $\cA \subseteq \cB$ such that $\ns{U} = \NSCup \cA$.
\end{definition}

\begin{notation}
\label{not:finite_intersections}
Given a family $\cS$ of neutrosophic sets over a common universe set $\U$,
let us denote by $\cB(\cS)$ (or, equivalently, by $\cS^*$) the \df{family of all its finite neutrosophic intersections},
with the addition of the absolute neutrosophic set, i.e:
$$\cB(\cS) = \left\{ \NSCap_{i=1}^n \ns{S}_i : \, n \in \NN, \ns{S}_i \in \cS \right\} \NScup \left\{ \NSabsoluteset \right\} .$$
\end{notation}

\begin{proposition}
\label{pro:finite_intersections_of_two_families}
If $\cS_1$ and $\cS_2$ are two families of neutrosophic sets on the same universe set $\U$
such that $\cS_1 \subseteq \cS_2$
then the corresponding family of all their finite neutrosophic intersections satisfy $\cB(\cS_1) \subseteq \cB(\cS_2)$.
\end{proposition}

\begin{proposition}
\label{pro:neutrosophic_topology_generated_by_a_subbase}
The neutrosophic topology $\Tau(\cS)$ generated by a family of neutrosophic sets $\cS$ over the same universe set $\U$
coincides with the topology $\Tau\left( \cB(\cS) \right)$ generated by the family
of all finite neutrosophic intersections $\cB(\cS)$ of $\cS$, i.e., we have:
$$\Tau(\cS) = \Tau\left( \cB(\cS) \right)$$
and the family $\cB(\cS)$ of all finite neutrosophic intersections forms a neutrosophic basis for $\Tau(\cS)$.
\end{proposition}

\begin{definition}
\label{def:neutrosophic_subbasis}
In the situation described by Proposition \ref{pro:neutrosophic_topology_generated_by_a_subbase},
the $\cS$ family is said to be a \df{neutrosophic sub-basis} (or a \df{neutrosophic pre-base}) of the $\Tau(\cS)$ topology.
\end{definition}


\subsection{The PYNS Framework}
The Python Neutrosophic Sets (PYNS) framework, developed by Nordo et al. \cite{nordo2024}, supports essential components for neutrosophic modeling.
It is composed of three primary classes: \cpy{NSuniverse}, \cpy{NSset} and \cpy{NSmapping}
designed to manage neutrosophic universes, individual neutrosophic sets and functions between neutrosophic sets, respectively.

The PYNS framework includes three primary classes that enable the representation and manipulation of neutrosophic structures in Python.
These are:
\begin{itemize}
    \item \texttt{NSuniverse} that defines a universe set over which neutrosophic sets are established. It supports initialization with various formats, such as lists, tuples, and strings. Key methods include:
    \begin{itemize}
        \item \texttt{get()} - Returns the elements of the universe set.
        \item \texttt{cardinality()} - Outputs the total count of elements.
        \item \texttt{isSubset(unv)} - Checks if the current universe set is a subset of another specified universe set.
    \end{itemize}

    \item \texttt{NSset} which represents individual neutrosophic sets, where each element in the universe has associated degrees of membership, indeterminacy, and non-membership. Primary methods include:
    \begin{itemize}
        \item \texttt{NSunion(nset)} - Computes the union of the current neutrosophic set with another specified set.
        \item \texttt{NSintersection(nset)} - Finds the intersection of the current set with another.
        \item \texttt{NScomplement()} - Generates the complement of the neutrosophic set, flipping the membership, indeterminacy, and non-membership values.
        \item \texttt{NSdifference(nset)} - Calculates the neutrosophic difference between the current set and another.
        \item \texttt{isNSdisjoint(nset)} - Verifies if the current set is disjoint from the given neutrosophic set.
    \end{itemize}

    \item \texttt{NSmapping} that is designed to manage mappings between two neutrosophic universe sets, each instance of this class holds a dictionary to represent the mapping relations. It includes methods such as:
    \begin{itemize}
        \item \texttt{getDomain()} - Retrieves the domain of the mapping as an \texttt{NSuniverse} object.
        \item \texttt{getCodomain()} - Returns the codomain as an \texttt{NSuniverse} object.
        \item \texttt{getMap()} - Provides access to the dictionary of element-value pairs defining the mapping.
        \item \texttt{setValue(u, v)} - Assigns a specific codomain element \(v\) to a domain element \(u\).
        \item \texttt{getValue(u)} - Returns the mapping value of a given element \(u\) from the domain.
        \item \texttt{getFibre(v)} - Identifies all elements in the domain that map to a specified codomain value \(v\), effectively finding the fibre of \(v\).
        \item \texttt{NSimage(nset)} - Computes the neutrosophic image of a given neutrosophic set by the mapping.
        \item \texttt{NScounterimage(nset)} - Determines the neutrosophic counterimage for the specified neutrosophic set.
    \end{itemize}
\end{itemize}

In summary, the \texttt{PYNS} framework, as developed by Nordo et al. \cite{nordo2024},
provides a robust foundation for neutrosophic modeling in Python, structured around
the core classes \texttt{NSuniverse}, \texttt{NSset}, and \texttt{NSmapping}.
With approximately $1500$ lines of code, \texttt{PYNS} supports intuitive operations
on neutrosophic universes, individual sets, and mappings.
Key functionalities include fundamental operations such as neutrosophic union, intersection,
difference, and complement, as well as the computation of images and counterimages through mappings,
all of which contribute to a flexible toolset for handling neutrosophic structures
across various types of universes.
The modular design of \texttt{PYNS} not only supports interactive use,
ideal for experimentation and the exploration of neutrosophic properties,
but also facilitates integration into larger Python projects due to its adaptable architecture
and comprehensive, example-rich documentation.
This ease of use and flexibility makes the \texttt{PYNS} framework a valuable open source tool
for exploring the properties of neutrosophic sets,
with applications spanning diverse areas of research that benefit from neutrosophic set theory.

Despite its comprehensive capabilities, \texttt{PYNS} currently lacks a mechanism for systematically managing collections of neutrosophic sets that could serve as bases or sub-bases in topological constructions. The ability to handle such collections is crucial for building and verifying topological structures in neutrosophic settings, which remain essential for advanced applications and theoretical explorations in neutrosophic topology.
This limitation is addressed by extending the PYNS framework with the \texttt{NSfamily} class,
presented in this paper which provides the necessary functionality to manage and interact
with families of neutrosophic sets and neutrosophic topologies.


\section{The \texttt{NSfamily} Class}
The \texttt{NSfamily} class enhances the \cpy{PYNS} framework by enabling the construction and management
of neutrosophic families, which serve as bases or sub-bases for topology generation.
Inspired by Ozturk's theoretical work on neutrosophic topology \cite{ozturk2020} and Salama et al.'s object-oriented design principles \cite{salama2014b, salama2014c}, \texttt{NSfamily} provides a unified interface for creating topologies
and verifying closure properties within neutrosophic environments.

By adding \texttt{NSfamily} to the other three classes \texttt{NSuniverse}, \texttt{NSset} and \texttt{NSmapping}
above described, the framework PYNS now supports topological operations on families of neutrosophic sets
within the same universe, bridging a key gap in the framework and enabling systematic topology construction.
\\
The UML diagram below illustrates the relationships among the classes, with dashed arrows indicating "uses" relationships.


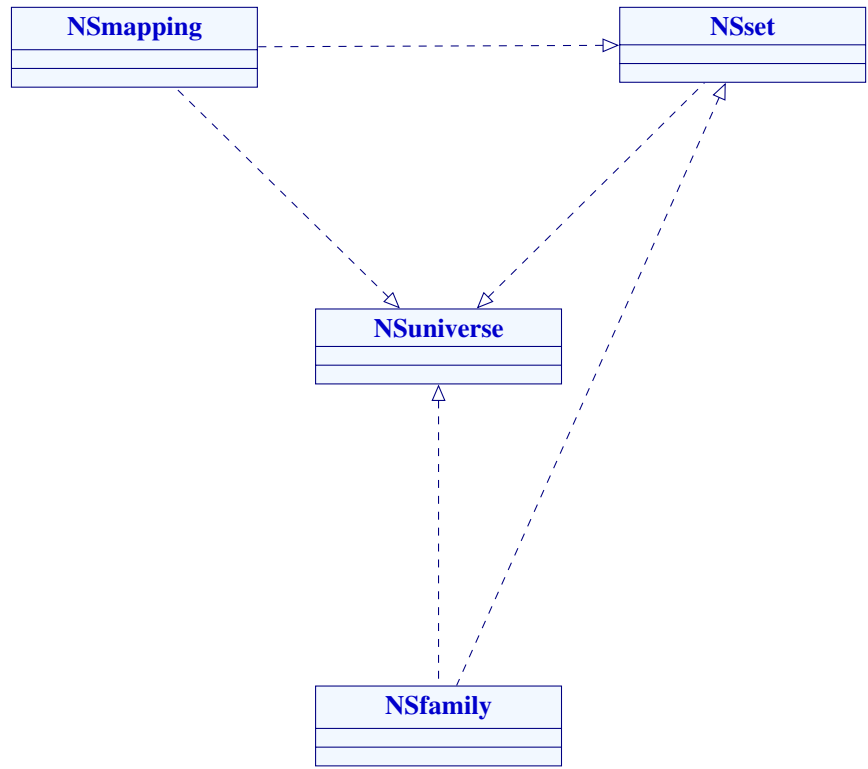
\vspace{2mm}
\begin{center}
\begin{tikzpicture}

  \begin{class}[text width=3cm]{NSuniverse}{4,5} 
  \end{class}

  \begin{class}[text width=3cm]{NSset}{8,9}   
    \implement{NSuniverse}
  \end{class}

  \begin{class}[text width=3cm]{NSmapping}{0,9}  
    \implement{NSuniverse}
    \implement{NSset}
  \end{class}

  \begin{class}[text width=3cm]{NSfamily}{4,0}
    \implement{NSuniverse}
    \implement{NSset}
  \end{class}

\end{tikzpicture}
\end{center}
\vskip 5mm


The \texttt{NSfamily} class manages families of neutrosophic sets stored as lists of objects of the class \texttt{NSset}
sharing the same universe set and can be shortly described
by means of its main properties and methods in the following UML class diagram.

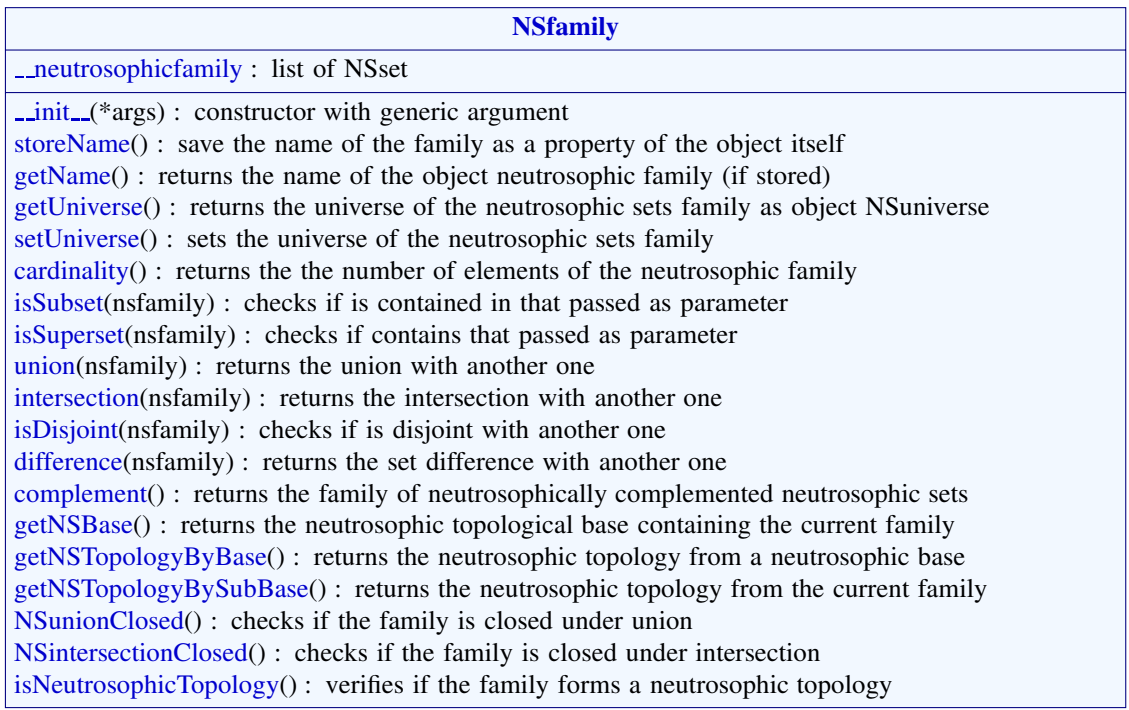
\begin{figure}[H]
\begin{center}
\begin{tikzpicture}
  \begin{class}[text width=14.5cm]{NSfamily}{0,0}
    \attribute{{\color{blu}\dunder neutrosophicfamily} : list of NSset}
    \operation{{\color{blu}\dunder init\dunder}(*args) : constructor with generic argument}
    \operation{{\color{blu}storeName}() : save the name of the family as a property of the object itself}
    \operation{{\color{blu}getName}() : returns the name of the object neutrosophic family (if stored)}
    \operation{{\color{blu}getUniverse}() : returns the universe of the neutrosophic sets family as object NSuniverse}
    \operation{{\color{blu}setUniverse}() : sets the universe of the neutrosophic sets family}
    \operation{{\color{blu}cardinality}() : returns the the number of elements of the neutrosophic family}
    \operation{{\color{blu}isSubset}(nsfamily) : checks if is contained in that passed as parameter}
    \operation{{\color{blu}isSuperset}(nsfamily) : checks if contains that passed as parameter}
    \operation{{\color{blu}union}(nsfamily) : returns the union with another one}
    \operation{{\color{blu}intersection}(nsfamily) : returns the intersection with another one}
    \operation{{\color{blu}isDisjoint}(nsfamily) : checks if is disjoint with another one}
    \operation{{\color{blu}difference}(nsfamily) : returns the set difference with another one}
    \operation{{\color{blu}complement}() : returns the family of neutrosophically complemented neutrosophic sets}
    \operation{{\color{blu}getNSBase}() : returns the neutrosophic topological base containing the current family}
    \operation{{\color{blu}getNSTopologyByBase}() : returns the neutrosophic topology from a neutrosophic base}
    \operation{{\color{blu}getNSTopologyBySubBase}() : returns the neutrosophic topology from the current family}
    \operation{{\color{blu}NSunionClosed}() : checks if the family is closed under union}
    \operation{{\color{blu}NSintersectionClosed}() : checks if the family is closed under intersection}
    \operation{{\color{blu}isNeutrosophicTopology}() : verifies if the family forms a neutrosophic topology}
 \end{class}
\end{tikzpicture}
\end{center}
\label{dia:nsfamily}
\end{figure}

The \texttt{NSfamily} constructor initializes a family of neutrosophic sets over a shared universe. The initialization process is detailed in the following algorithm:

\begin{algoritmo}
\caption{Constructor method of the class \cpy{NSfamily}}
\SetKwFunction{MyFunction}{$\dunder$init$\dunder$}
\SetKwProg{Fn}{Function}{:}{}
\Fn{\MyFunction{args}}{
    Initialize the private property \textit{neutrosophicfamily} as an empty list\;
    Get the \textit{length} of \textit{args}\;

    \If{\textit{length} $= 0$}{
        Set \textit{universe} to \texttt{None}
    }

    \ElseIf{\textit{length} $= 1$}{
        Let \textit{elem} be the first and only element in \textit{args}\;

        \If{\textit{elem} \textbf{is an instance of} \cpy{NSset}}{
            Add \textit{elem} to \textit{neutrosophicfamily} and set \textit{universe} to the universe of \textit{elem}\;
        }
        \ElseIf{\textit{elem} \textbf{is a list or tuple}}{
            \ForEach{\textit{e} in \textit{elem}}{
                \If{\textit{e} is not in \textit{neutrosophicfamily}}{
                    Set the name of \textit{e} and add it to \textit{neutrosophicfamily}\;
                }
            }
            \If{\textit{neutrosophicfamily} is not empty}{
                Set \textit{universe} to the universe of the first element in \textit{neutrosophicfamily}\;
            }
            \Else{
                Set \textit{universe} to \texttt{None}\;
            }
        }
    }

    \ElseIf{\textit{length} $> 1$}{
        \ForEach{\textit{e} in \textit{args}}{
            \If{\textit{e} is not in \textit{neutrosophicfamily}}{
                Set the name of \textit{e}  and add it to \textit{neutrosophicfamily}\;
            }
        }
        Set \textit{universe} to the universe of the first element in \textit{args}\;
    }

    Store \textit{universe}, \textit{neutrosophicfamily} as private properties\;
    Set the private property \textit{name} to \texttt{None}\;
}
\end{algoritmo}

\vskip 4mm

The Python code corresponding to the constructor method of the \texttt{NSfamily} class is given below.

\begin{framed}
\begin{codpy}
from .ns_universe import NSuniverse
from .ns_set import NSset
from .ns_util import NSreplace,NSstringToDict,NSisExtDict,nameToBB,isBB
import inspect
from itertools import combinations
from functools import reduce
from time import time

class NSfamily:

    def $\dunder$init$\dunder$(self, *args):
        neutrosophicfamily = list()
        length = len(args)
        if length == 0:
            universe = None
        elif length == 1:
            elem = args[0]
            if type(elem) == NSset:
                neutrosophicfamily = [elem]
                universe = elem.getUniverse()
            elif type(elem) in [list ,tuple]:
                for e in elem:
                    if e not in neutrosophicfamily:
                        e.setName(e.getName())
                        neutrosophicfamily.append(e)
                if len(neutrosophicfamily) > 0:
                    universe = neutrosophicfamily[0].getUniverse()
                else:
                    universe = None
        elif length > 1:
            for e in args:
                if e not in neutrosophicfamily:
                    e.setName(e.getName())
                    neutrosophicfamily.append(e)
            universe = args[0].getUniverse()
        self.$\dunder$universe = universe
        self.$\dunder$neutrosophicfamily = neutrosophicfamily
        self.$\dunder$name = None
\end{codpy}
\end{framed}


The \texttt{NSfamily} constructor offers flexible initialization to handle various cases
in creating neutrosophic set families, managing arguments through the \texttt{*args} parameter.
This enables creation from individual sets, collections, or copies of existing \texttt{NSfamily} objects:
\begin{itemize}
    \item empty family: \texttt{NSfamily()} creates an empty family with no universe, allowing elements to be added incrementally.

    \item single neutrosophic set: \texttt{NSfamily(ns\_set)} initializes a family with one set, inheriting its universe.

    \item list or tuple of sets: \texttt{NSfamily([ns\_set1, ns\_set2, \ldots])} or or \texttt{NSfamily((ns\_set1, ns\_set2, \ldots))}
    accepts a collection, avoiding duplicates and setting the universe based on the first set.

    \item multiple sets as arguments: \texttt{NSfamily(ns\_set1, ns\_set2, \ldots)} creates a family
    from the specified sets, with the universe derived from the first neutrosophic set,
    ensuring no duplicate entries.

    \item existing \texttt{NSfamily} object: \texttt{NSfamily(existing\_family)} copies an existing family, creating an independent object with the same universe.
\end{itemize}
This versatility allows the \texttt{NSfamily} class to support a broad range of initial configurations,
thereby facilitating the manipulation and management of neutrosophic set families.


In order to conveniently display families of neutrosophic sets in both simplified and tabular text formats while providing a comprehensive representation, the special methods \cpy{$\dunder$str$\dunder$()} and \cpy{$\dunder$format$\dunder$()} have been implemented. These methods allow customization of the output format, especially when using the \cpy{print} function with f-strings
for displaying \cpy{NSfamily} objects.
The \cpy{$\dunder$format$\dunder$()} method supports the following format specifiers:
\begin{description}[style=unboxed,labelindent=1em]
    \item[\textbf{s}] : simple format (default), offering a straightforward textual representation
    \item[\textbf{t}] : tabular format, presenting data in a structured and neatly aligned manner
    \item[\textbf{l}] : includes a label with the \cpy{NSfamily} object's name, if available
    \item[\textbf{x}] : extended format, providing additional spacing in the first column to accommodate longer names or labels.
\end{description}

\begin{framed}
\begin{codpy}
def $\dunder$str$\dunder$(self, tabularFormat=False, label=False, extended=False):
    labelname = ""
    if label and self.$\dunder$name is not None:
        labelname = f"{self.$\dunder$name} = "
    if not self.$\dunder$neutrosophicfamily:
        return labelname + "\u2205"
    indentation = " "*(len(labelname) + 2)
    if not tabularFormat:
        if extended:
            items = [a.$\dunder$str$\dunder$(tabularFormat, True, extended)
                     for a in self.$\dunder$neutrosophicfamily]
            formatted_lines=[", ".join(items[i:i + 2]) for i in range(0,len(items),2)]
            return labelname+"{ "+f",\n{indentation}".join(formatted_lines)+" }\n"
        else:
            return labelname + "{ " + ", ".join(
                [a.$\dunder$str$\dunder$(tabularFormat, True, extended)
                 for a in self.$\dunder$neutrosophicfamily]) + " }\n"
    res = labelname + "{ " + ", ".join(
        [a.$\dunder$str$\dunder$(tabularFormat, True, extended)
         for a in self.$\dunder$neutrosophicfamily]) + " }\n"
    return res

def $\dunder$format$\dunder$(self, spec):
    label = "l" in spec
    extended = "x" in spec
    if "t" in spec:
        result = self.$\dunder$str$\dunder$(tabularFormat=True, label=label, extended=extended)
    else:
        result = self.$\dunder$str$\dunder$(tabularFormat=False, label=label, extended=extended)
    return result
\end{codpy}
\end{framed}

Some basic methods facilitate the management of neutrosophic family properties.
\\
The \cpy{storeName()} method captures the object's name from the local context,
while \cpy{getName()} retrieves it if stored. \cpy{getUniverse()} returns the associated \cpy{NSuniverse} of the neutrosophic family, and \cpy{setUniverse()} sets it, ensuring it is a valid \cpy{NSuniverse} instance.
Finally, \cpy{cardinality()} provides the count of neutrosophic sets within the family.
Their code is given below.

\begin{framed}
\begin{codpy}
def storeName(self):
    frame = inspect.currentframe().f_back
    local_vars = frame.f_locals
    var_name = next((nm for nm, vl in local_vars.items()
                        if vl is self), None)
    self.__name = var_name

def getName(self):
    return self.__name

def getUniverse(self):
    return self.__universe

def setUniverse(self, universe):
    if type(universe) != NSuniverse:
        raise ValueError("The parameter's type must be NSuniverse.")
    self.__universe = universe

def cardinality(self):
    return len(self.__neutrosophicfamily)
\end{codpy}
\end{framed}


Let us illustrate what just said with an example of code
executed interactively in the Python console.

\begin{codpyconsole}
>>> from NS.pyns.ns_universe import NSuniverse
>>> from NS.pyns.ns_set import NSset
>>> from NS.pyns.ns_family import NSfamily
>>> U = NSuniverse("a,b,c")
>>> A1 = NSset(U, "(0.4,0.4,0.3), (0.1,0.1,0.1), (0.2,0.2,0.2)")
>>> A2 = NSset(U, "(0.1,0.2,0.9), (0.9,0.1,0.3), (0.5,0.3,0.4)")
>>> A1.storeName()
>>> A2.storeName()
>>> E = NSset.EMPTY(U)
>>> L = NSfamily(E, A1, A2)
>>> L.storeName()
>>> print(f"family {L:lx}")
family L = { A1 = < a/(0.4,0.4,0.3), b/(0.1,0.1,0.1), c/(0.2,0.2,0.2) >,
      A2 = < a/(0.1,0.2,0.9), b/(0.9,0.1,0.3), c/(0.5,0.3,0.4) > }
\end{codpyconsole}
\vskip 2mm


The following methods provide essential operations for managing and comparing neutrosophic families.
These methods allow for determining subset and superset relationships, computing unions and intersections of families, checking for disjoint sets, and calculating differences and complements.

The \cpy{isSubset()} method verifies if the current neutrosophic family is contained within another, while the \cpy{__le__()} method overloads the \cpy{<=} operator to provide a more intuitive comparison.

\begin{framed}
\begin{codpy}
def isSubset(self, nsfamily):
    if type(nsfamily) != NSfamily:
        raise ValueError("the parameter is not a neutrosophic family")
    if self.getUniverse() != nsfamily.getUniverse():
        raise ValueError("the two neutrosophic families cannot be defined on different universe sets")
    else:
        result = True
        for e in self.__neutrosophicfamily:
            if e not in nsfamily.__neutrosophicfamily:
                result = False
                break
        return result

def __le__(self, nsfamily):
    if type(nsfamily) != NSfamily:
        raise ValueError("the second argument is not a neutrosophic family")
    return self.isSubset(nsfamily)
\end{codpy}
\end{framed}

Similarly, the \cpy{isSuperset()} method checks if the current neutrosophic family contains another, with the \cpy{__ge__()} method overloading the \cpy{>=} operator to simplify this comparison.

\begin{framed}
\begin{codpy}
def isSuperset(self, nsfamily):
    if type(nsfamily) != NSfamily:
        raise ValueError("the parameter is not a neutrosophic family")
    if self.getUniverse() != nsfamily.getUniverse():
        raise ValueError("the two neutrosophic families cannot be defined on different universe sets")
    return nsfamily.isSubset(self)

def __ge__(self, nsfamily):
    if type(nsfamily) != NSfamily:
        raise ValueError("the second argument is not a neutrosophic family")
    return self.isSuperset(nsfamily)
\end{codpy}
\end{framed}


Equality and inequality between families are determined by overloading
the comparison operators \cpy{==} and \cpy{!=} through the special methods \cpy{__eq__()} and \cpy{__ne__()}.

\begin{framed}
\begin{codpy}
def __eq__(self, nsfamily):
    if type(nsfamily) != NSfamily:
        raise ValueError("the second argument is not a neutrosophic family")
    if self.getUniverse() != nsfamily.getUniverse():
        raise ValueError("the two neutrosophic families cannot be defined on different universe sets")
    equal = self.isSubset(nsfamily) and nsfamily.isSubset(self)
    return equal

def __ne__(self, nsfamily):
    if type(nsfamily) != NSfamily:
        raise ValueError("the second argument is not a neutrosophic family")
    if self.getUniverse() != nsfamily.getUniverse():
        raise ValueError("the two neutrosophic families cannot be defined on different universe sets")
    different = not (self == nsfamily)
    return different
\end{codpy}
\end{framed}


The \cpy{union()} method, along with its operator \cpy{+} overloaded by \cpy{__add__()}, combines two neutrosophic families into one, eliminating duplicates and ensuring they share the same universe.

\begin{framed}
\begin{codpy}
def union(self, nsfamily):
    if not isinstance(nsfamily, NSfamily):
        raise ValueError("The parameter is not a neutrosophic family.")
    if self.getUniverse() != nsfamily.getUniverse():
        raise ValueError("The two neutrosophic families cannot be defined
            on different universes.")
    combined_sets = self.__neutrosophicfamily +
            [s for s in nsfamily.__neutrosophicfamily
                    if s not in self.__neutrosophicfamily]
    return NSfamily(combined_sets)

def __add__(self, nsfamily):
    if not isinstance(nsfamily, NSfamily):
        raise ValueError("the second argument is not a neutrosophic family")
    return self.union(nsfamily)
\end{codpy}
\end{framed}

The \cpy{intersection()} method identifies common elements between two families, facilitated by the \cpy{\&} operator via \cpy{__and__()}, while \cpy{isDisjoint()} checks for disjointness.

\begin{framed}
\begin{codpy}
def intersection(self, nsfamily):
    if not isinstance(nsfamily, NSfamily):
        raise ValueError("The parameter is not a neutrosophic family.")
    if self.getUniverse() != nsfamily.getUniverse():
        raise ValueError("The two neutrosophic families cannot be defined on different universes.")
    common_sets = [s for s in self.__neutrosophicfamily
                if s in nsfamily.__neutrosophicfamily]
    return NSfamily(common_sets)

def __and__(self, nsfamily):
    if not isinstance(nsfamily, NSfamily):
        raise ValueError("the second argument is not a neutrosophic family")
    return self.intersection(nsfamily)

def isDisjoint(self, nsfamily):
    nsemptyfamily = NSfamily()
    nsemptyfamily.setUniverse(self.getUniverse())
    intersez = self.intersection(nsfamily)
    disjoint = intersez.cardinality() == 0
    return disjoint
\end{codpy}
\end{framed}

The \cpy{difference()} method and its \cpy{-} operator via \cpy{__sub__()} return the elements unique to the first family, and \cpy{complement()} along with \cpy{\~} (overloaded by the special method
\cpy{__invert__()}) yields the family consisting of the neutrosophic complements of the neutrosophic sets of the original family.

the family of all complement of the current family.

\begin{framed}
\begin{codpy}
def difference(self, nsfamily):
    if not isinstance(nsfamily, NSfamily):
        raise ValueError("The parameter is not a neutrosophic family.")
    if self.getUniverse() != nsfamily.getUniverse():
        raise ValueError("The two neutrosophic families cannot be defined on different universes.")
    unique_sets = [s for s in self.__neutrosophicfamily
                    if s not in nsfamily.__neutrosophicfamily]
    return NSfamily(unique_sets)

def __sub__(self, nsfamily):
    if not isinstance(nsfamily, NSfamily):
        raise ValueError("The second argument is not a neutrosophic family")
    return self.difference(nsfamily)

def complement(self):
    complementary_sets = [ns_set.NScomplement()
                            for ns_set in self.__neutrosophicfamily]
    return NSfamily(complementary_sets)

def __invert__(self):
    return self.complement()
\end{codpy}
\end{framed}


\begin{codpyconsole}
>>> from pyns.ns_universe import NSuniverse
>>> from pyns.ns_set import NSset
>>> from pyns.ns_family import NSfamily
>>> U = NSuniverse("a,b,c")
>>> A1 = NSset(U, "(0.4,0.4,0.3), (0.1,0.1,0.1), (0.2,0.2,0.2)")
>>> A2 = NSset(U, "(0.1,0.2,0.9), (0.9,0.1,0.3), (0.5,0.3,0.4)")
>>> A3 = NSset(U, "(0.7,0.3,0.1), (0.8,0.4,0.0), (0.1,0.1,0.9)")
>>> A4 = NSset(U, "(0.2,0.2,0.8), (0.6,0.6,0.3), (0.5,0.4,0.5)")
>>> A1.storeName()
>>> A2.storeName()
>>> A3.storeName()
>>> A4.storeName()
>>> L1 = NSfamily(A1, A2, A3)
>>> L2 = NSfamily(A3, A4)
>>> L1.storeName()
>>> L2.storeName()
>>> L3 = L1 & L2
>>> L3.storeName()
>>> print(f"{L3:lx}")
L3 = { A3 = < a/(0.7,0.3,0.1), b/(0.8,0.4,0), c/(0.1,0.1,0.9) > }
>>> print(L3 <= L1)
True
>>> L4 = L1 + L2
>>> L4.storeName()
>>> print(f"{L4:lx}")
L4 = { A1 = < a/(0.4,0.4,0.3), b/(0.1,0.1,0.1), c/(0.2,0.2,0.2) >,
       A2 = < a/(0.1,0.2,0.9), b/(0.9,0.1,0.3), c/(0.5,0.3,0.4) >,
       A3 = < a/(0.7,0.3,0.1), b/(0.8,0.4,0), c/(0.1,0.1,0.9) >,
       A4 = < a/(0.2,0.2,0.8), b/(0.6,0.6,0.3), c/(0.5,0.4,0.5) > }
>>> L5 = $\tld$L4
>>> L5.storeName()
>>> print(f"{L5:lx}")
L5 = { $\tld$A1 = < a/(0.3,0.6,0.4), b/(0.1,0.9,0.1), c/(0.2,0.8,0.2) >,
       $\tld$A2 = < a/(0.9,0.8,0.1), b/(0.3,0.9,0.9), c/(0.4,0.7,0.5) >,
       $\tld$A3 = < a/(0.1,0.7,0.7), b/(0,0.6,0.8), c/(0.9,0.9,0.1) >,
       $\tld$A4 = < a/(0.8,0.8,0.2), b/(0.3,0.4,0.6), c/(0.5,0.6,0.5) > }
>>> L4.isDisjoint(L5)
True
\end{codpyconsole}


From Proposition \ref{pro:neutrosophic_topology_generated_by_a_subbase}
and Definition \ref{def:neutrosophic_subbasis}
we know that any family of neutrosophic sets over a same universe set
is a neutrosophic sub-basis for the neutrosophic topology $\Tau(\cS)$ generated by $\cS$
which has the family $\cB(\cS)$ of all finite neutrosophic intersections of $\cS$
as neutrosophic basis.
\\
The method \texttt{getNSBase} returns a neutrosophic topological basis from
the current neutrosophic family regarded as a sub-basis
by generating all its possible neutrosophic intersections.

\begin{algoritmo}
\caption{Method \texttt{getNSBase}}
\SetKwFunction{MyFunction}{getNSBase}
\SetKwProg{Fn}{Function}{:}{}
\Fn{\MyFunction{}}{
    Let \textit{subbase} be \texttt{self.\_\_neutrosophicfamily}\;
    Initialize an empty list \textit{base}\;

    \For{$i \leftarrow 1$ \KwTo \texttt{len(subbase)}}{
        \ForEach{\textit{combination} in \texttt{combinations(subbase, i)}}{
            Compute \textit{intersection} by applying \texttt{NSintersection} on all elements in \textit{combination}\;

            \If{\textit{intersection\_name} is not empty}{
                Set the name of \textit{intersection} to \textit{intersection\_name}\;
            }

            \ForEach{\textit{s} in \textit{subbase}}{
                \If{\textit{intersection} equals \textit{s}}{
                    Set the name of \textit{intersection} to name of \textit{s}\;
                    \textbf{break}\;
                }
            }
            Add \textit{intersection} to \textit{base}\;
        }
    }
    Convert \textit{base} to an \texttt{NSfamily} object\;
    Set the universe of \textit{base} to \texttt{self.getUniverse()}\;
    \Return{\textit{base}}\;
}
\end{algoritmo}
\vskip 5mm

The Python code corresponding to the method \texttt{getNSBase} is given below.

\begin{framed}
\begin{codpy}
def getNSBase(self):
    subbase = self.$\dunder$neutrosophicfamily
    base = list()
    for i in range(1, len(subbase) + 1):
        for combin in combinations(subbase, i):
            intersez = reduce(lambda x, y: x.NSintersection(y), combin)
            names = [s.getName() for s in combin if s.getName()]
            intersez_name = " $\NScap$ ".join(names) if names else None
            if intersez_name:
                intersez.setName(intersez_name)
            for s in subbase:
                if intersez == s:
                    intersez.setName(s.getName())
                    break
            base.append(intersez)
    base = NSfamily(base)
    base.setUniverse(self.getUniverse())
    return base
\end{codpy}
\end{framed}
\vskip 4mm

As an example, we show how the method described above
can be used in the interactive mode by means of the Python console:

\begin{codpyconsole}
>>> from pyns.ns_universe import NSuniverse
>>> from pyns.ns_set import NSset
>>> from pyns.ns_family import NSfamily
>>> U = NSuniverse("a,b,c")
>>> A1 = NSset(U, "(0.4,0.4,0.3), (0.1,0.1,0.1), (0.2,0.2,0.2)")
>>> A2 = NSset(U, "(0.1,0.2,0.9), (0.9,0.1,0.3), (0.5,0.3,0.4)")
>>> A1.storeName()
>>> A2.storeName()
>>> L = NSfamily(A1, A2)
>>> L.storeName()
>>> print(f"the neutrosophic family is {L:lt}")
 the neutrosophic family is L = {
 A1      |   membership   |  indeterminacy | non-membership |
-------------------------------------------------------------
 a       | 0.4            | 0.4            | 0.3            |
 b       | 0.1            | 0.1            | 0.1            |
 c       | 0.2            | 0.2            | 0.2            |
-------------------------------------------------------------
,
 A2      |   membership   |  indeterminacy | non-membership |
-------------------------------------------------------------
 a       | 0.1            | 0.2            | 0.9            |
 b       | 0.9            | 0.1            | 0.3            |
 c       | 0.5            | 0.3            | 0.4            |
-------------------------------------------------------------
>>> B = L.getNSBase()
>>> B.storeName()
>>> print(f"basis is {B:lx}")
basis is B = { A1 = < a/(0.4,0.4,0.3), b/(0.1,0.1,0.1), c/(0.2,0.2,0.2) >,
      A2 = < a/(0.1,0.2,0.9), b/(0.9,0.1,0.3), c/(0.5,0.3,0.4) >,
      A1 $\NScap$ A2 = < a/(0.1,0.2,0.9), b/(0.1,0.1,0.3), c/(0.2,0.2,0.4) > }
\end{codpyconsole}
\vskip 2mm


From Proposition \ref{pro:neutrosophic_topology_generated_by_a_subbase}
and Definition \ref{def:neutrosophicbase},
we know that a base for a topology is a collection of sets from which
every open set can be represented as a union of elements of the base.
Specifically, in the context of neutrosophic topology, a neutrosophic base
provides the foundation for generating all possible neutrosophic unions, which form the topology.
\\
The method \texttt{getNSTopologyByBase} returns the neutrosophic topology
derived from the current base. This is accomplished by generating all possible
finite neutrosophic unions of the base sets, thus constructing the full topology.
\vskip 2mm

\begin{algoritmo}
\caption{Method \texttt{getNSTopologyByBase}}
\SetKwFunction{MyFunction}{getNSTopologyByBase}
\SetKwProg{Fn}{Function}{:}{}
\Fn{\MyFunction{}}{
    Let \textit{base} be \texttt{self.\_\_neutrosophicfamily}\;
    Initialize an empty list \textit{topology}\;
    Create an \textit{empty} neutrosophic set and add it to \textit{topology}\;

    \For{$i \leftarrow 1$ \KwTo \texttt{len(base)}}{
        \ForEach{\textit{combination} in \texttt{combinations(base, i)}}{
            Compute \textit{union} by applying \texttt{NSunion} on all elements in \textit{combination}\;
            Construct the \textit{union\_name} based on the combination's names\;
            Assign \textit{union\_name} to \textit{union}, if applicable\;
            \If{\textit{union} matches any \textit{base} element}{
                Update \textit{union}'s name accordingly\;
            }
            Add \textit{union} to \textit{topology}\;
        }
    }
    Create an \textit{absolute} neutrosophic set from the universe and add it to \textit{topology}\;
    Convert \textit{topology} to an \texttt{NSfamily} object\;
    Set the universe of \textit{topology} to \texttt{self.getUniverse()}\;
    \Return{\textit{topology}}\;
}
\end{algoritmo}
\vskip 3mm

The Python code corresponding to the method \texttt{getNSTopologyByBase} is given below.

\begin{framed}
\begin{codpy}
def getNSTopologyByBase(self):
    base = self.$\dunder$neutrosophicfamily
    topology = list()
    universe = self.$\dunder$universe
    empty = NSset.EMPTY(universe)
    empty.setName("\u2205\u0303")
    topology.append(empty)
    for i in range(1, len(base) + 1):
        for combin in combinations(base, i):
            union = reduce(lambda x, y: x.NSunion(y), combin)
            names = [s.getName() for s in combin if s.getName()]
            union_name = " $\NScup$ ".join(f"({name})" if "$\NScap$" in name else name
                                for name in names)
            if union_name:
                union.setName(union_name)
            for b in base:
                if union == b:
                    union.setName(b.getName())
                    break
            topology.append(union)
    absolute = NSset.ABSOLUTE(universe)
    universe_name = universe.getName()
    absolute.setName(nameToBB(universe_name) if universe_name
                                and not isBB(universe_name) else universe_name)
    topology.append(absolute)
    topology = NSfamily(topology)
    topology.setUniverse(self.getUniverse())
    return topology
\end{codpy}
\end{framed}


The \texttt{getNSTopologyBySubBase} method constructs the neutrosophic topology derived from a subbase by computing all possible neutrosophic unions of intersections. It first generates the neutrosophic base using the \texttt{getNSBase()} method and then forms the complete topology by applying \texttt{getNSTopologyByBase()} on this base. This process ensures the creation of a full neutrosophic topology from any given family of neutrosophic sets.
The Python code corresponding to this method is shown below.

\begin{framed}
\begin{codpy}
def getNSTopologyBySubBase(self):
    nsbase = self.getNSBase()
    nstopology = nsbase.getNSTopologyByBase()
    return nstopology
\end{codpy}
\end{framed}


The \texttt{isNeutrosophicTopology} method verifies whether a neutrosophic family satisfies the axioms of a neutrosophic topology. It ensures the presence of both the empty set and the universal set in the family, and checks closure properties under union and intersection operations.

\begin{framed}
\begin{codpy}
def isNeutrosophicTopology(self):
    family = self.__neutrosophicfamily
    universe = self.__universe
    empty = NSset.EMPTY(universe)
    if empty not in family:
        return False
    absolute = NSset.ABSOLUTE(universe)
    if absolute not in family:
        return False
    if not self.NSunionClosed():
        return False
    if not self.NSintersectionClosed():
        return False
    return True
\end{codpy}
\end{framed}

\noindent
The \texttt{isNeutrosophicTopology} method relies on \texttt{NSunionClosed}
and \texttt{NSintersectionClosed}, both of which delegate the closure verification
to the private method \texttt{$\dunder$checkClosure}.
This method efficiently manages the verification process by iterating over
all possible combinations of sets and applying the specified operation (union or intersection)
provided as a lambda function. By centralizing the core logic
of closure checks, \texttt{$\dunder$checkClosure} ensures consistency and efficiency
in determining whether the neutrosophic family satisfies the closure properties
essential for a neutrosophic topology.
\vskip 2mm

\begin{framed}
\begin{codpy}
def __checkClosure(self, operation, operation_name):
    family = self.$\dunder$neutrosophicfamily
    l = len(family)
    for i in range(2, l + 1):
        for combin in combinations(family, i):
            result = reduce(operation, combin)
            if result not in family:
                return False
    return True


    def NSunionClosed(self):
        return self.__checkClosure(lambda x, y: x.NSunion(y), "union")

    def NSintersectionClosed(self:
        return self.__checkClosure(lambda x, y: x.NSintersection(y), "intersection")
\end{codpy}
\end{framed}

The following code example executed interactively in the Python console
illustrates the use of the above described methods.

\begin{codpyconsole}
>>> from pyns.ns_universe import NSuniverse
>>> from pyns.ns_set import NSset
>>> from pyns.ns_family import NSfamily
>>> U = NSuniverse("1,2,3")
>>> B1 = NSset(U, "(0.2,0.4,0.3), (0.6,0.1,0.1), (0.4,0.6,0.3)")
>>> B2 = NSset(U, "(0.3,0.2,0.9), (0.6,0.5,0.3), (0.2,0.3,0.8)")
>>> B1.storeName()
>>> B2.storeName()
>>> S = NSfamily(B1, B2)
>>> T = S.getNSTopologyBySubBase()
>>> T.storeName()
>>> print(f"topology has cardinality {T.cardinality()} and is:\n {T:tlx}")
topology has cardinality 6 and is:
 T = {
 $\nsempty$            |   membership   |  indeterminacy | non-membership |
------------------------------------------------------------------
 1            | 0              | 0              | 1              |
 2            | 0              | 0              | 1              |
 3            | 0              | 0              | 1              |
------------------------------------------------------------------
,
 B1           |   membership   |  indeterminacy | non-membership |
------------------------------------------------------------------
 1            | 0.2            | 0.4            | 0.3            |
 2            | 0.6            | 0.1            | 0.1            |
 3            | 0.4            | 0.6            | 0.3            |
------------------------------------------------------------------
,
 B2           |   membership   |  indeterminacy | non-membership |
------------------------------------------------------------------
 1            | 0.3            | 0.2            | 0.9            |
 2            | 0.6            | 0.5            | 0.3            |
 3            | 0.2            | 0.3            | 0.8            |
------------------------------------------------------------------
,
 B1 $\NScap$ B2      |   membership   |  indeterminacy | non-membership |
------------------------------------------------------------------
 1            | 0.2            | 0.2            | 0.9            |
 2            | 0.6            | 0.1            | 0.3            |
 3            | 0.2            | 0.3            | 0.8            |
------------------------------------------------------------------
,
 B1 $\NScup$ B2      |   membership   |  indeterminacy | non-membership |
------------------------------------------------------------------
 1            | 0.3            | 0.4            | 0.3            |
 2            | 0.6            | 0.5            | 0.1            |
 3            | 0.4            | 0.6            | 0.3            |
------------------------------------------------------------------
,
 $\nsabsol$            |   membership   |  indeterminacy | non-membership |
------------------------------------------------------------------
 1            | 1              | 1              | 0              |
 2            | 1              | 1              | 0              |
 3            | 1              | 1              | 0              |
------------------------------------------------------------------
 }
>>> T.isNeutrosophicTopology()
True
\end{codpyconsole}


\section{Conclusions}

In this paper, we introduced an extension to the Python Neutrosophic Sets (PYNS) framework
with the development of the \texttt{NSfamily} class.
This addition allows for the efficient construction, manipulation,
and analysis of neutrosophic families and their roles as bases and sub-bases
in neutrosophic topological spaces. Leveraging existing PYNS classes, \texttt{NSuniverse} and \texttt{NSset},
the new class provides robust methods for evaluating closure properties, unions, intersections, and complements,
which are essential operations in topological studies.

The theoretical exploration underscored the adaptability of neutrosophic sets
in addressing indeterminate and inconsistent information, a feature that differentiates them from classical set theories.
The practical implementations demonstrated through detailed examples reveal
the computational efficiency and flexibility of the PYNS framework,
particularly in generating complete neutrosophic topologies from defined families.

Future research will aim to broaden the applicability of neutrosophic topologies,
especially in decision-making scenarios characterized by high uncertainty and indeterminacy.
Further enhancements to the PYNS framework are also planned to optimize performance
and extend support for more complex neutrosophic operations.

This work lays a foundational framework for the integration of neutrosophic topological
methods into various scientific and engineering disciplines.
By providing a comprehensive toolset for the exploration and application of neutrosophic concepts,
it opens new avenues for research where traditional approaches may not suffice.

The source code for the PYNS framework, including the \texttt{NSfamily} class introduced in this paper,
is available at \url{https://github.com/giorgionordo/PythonNeutrosophicTopologies}.
This repository, licensed under GPL 3.0, serves as a valuable resource for researchers and practitioners,
enabling them to explore, utilize, and extend neutrosophic topologies in their projects.
\vskip 4mm


\noindent
\textbf{Acknowledgements}
This research was supported by Gruppo Nazionale per le Strutture Algebriche, Geometriche e le loro Applicazioni
(G.N.S.A.G.A.) of Istituto Nazionale di Alta Matematica (INdAM) "F.~Severi", Italy.


\end{document}